\documentclass[12pt]{article}
\usepackage{amssymb}
\usepackage{latexsym,bm}
\usepackage{graphicx}
\usepackage{amsmath}
\usepackage{mathrsfs}

\date{}
\begin{document}
\title{Graphs with many independent vertex cuts\footnote{E-mail addresses:
{\tt huyanan530@163.com}(Y.Hu),
{\tt zhan@math.ecnu.edu.cn}(X.Zhan), {\tt mathdzhang@163.com}(L.Zhang).}}
\author{\hskip -10mm Yanan Hu, Xingzhi Zhan and Leilei Zhang\thanks{Corresponding author.}\\
{\hskip -10mm \small Department of Mathematics, East China Normal University, Shanghai 200241, China}}\maketitle
\begin{abstract}
 The cycles are the only $2$-connected graphs in which any two nonadjacent vertices form a vertex cut. We generalize this fact by proving that
 for every integer $k\ge 3$ there exists a unique graph $G$ satisfying the following conditions: (1) $G$ is $k$-connected; (2) the independence
 number of $G$ is greater than $k;$ (3) any independent set of cardinality $k$ is a vertex cut of $G.$ The edge version of this result does not
 hold.  We also consider the problem when replacing independent sets by the periphery.
\end{abstract}



We consider finite simple graphs. For terminology and notations we follow the books [2,\,5]. It is known [4, p.46] that the cycles are the only $2$-connected graphs in which any two nonadjacent vertices form a vertex cut. We will generalize this fact and consider two related problems.

We denote by $V(G)$ the vertex set of a graph $G.$ The order of $G,$ denoted by $|G|,$ is the number of vertices of $G.$
For $S\subseteq V(G),$ the notation $G[S]$ denotes the subgraph of $G$ induced by $S.$
Let $K_{s,\, t}$ denote the complete bipartite graph whose partite sets have cardinality $s$ and $t,$ respectively.

{\bf Notation.} The notation $K_{s,s}-PM$ denotes the graph obtained from the balanced complete bipartite graph $K_{s,s}$
by deleting all the edges in a perfect matching of $K_{s,s}.$

Note that $K_{s,s}-PM$ is an $(s-1)$-connected $(s-1)$-regular graph, $K_{3,3}-PM$ is the $6$-cycle $C_6$ and  $K_{4,4}-PM$ is the cube $Q_3.$

{\bf Theorem 1.} {\it Let $k\ge 3$ be an integer. Then $K_{k+1,k+1}-PM$ is the  unique graph $G$ satisfying the following three
conditions: (1) $G$ is $k$-connected; (2) the independence  number of $G$ is greater than $k;$ (3) any independent set of cardinality $k$
is a vertex cut of $G.$}

{\bf Proof.} Let $G$ be a graph satisfying the three conditions in Theorem 1. We first assert that $G$ has order at least $2k+2.$
Let $S$ be an independent set of $G$ with cardinality $k+1.$ Since $G$ is $k$-connected, every vertex has degree at least $k.$
Let $T$ be the neighborhood of one vertex in $S.$ Then $|T|\ge k.$  Thus $|G|\ge |S|+|T|\ge 2k+1.$
If $|G|=2k+1,$ then $T$ would be the common neighborhood of all the vertices in $S.$ But now any $k$ vertices in $S$ do not form a vertex cut,
contradicting condition (3). This shows that $|G|\ge 2k+2.$

Choose an arbitrary but fixed independent set $A=\{x_1,x_2,\ldots,x_{k+1}\}$ of cardinality $k+1$ in $G.$ By condition (3), 
for every $i$ with $1\le i\le k+1,$ the graph $H_i\triangleq G-(A\setminus\{x_i\})$ is disconnected. Let $G_i$ denote the union of all the components of $H_i$ except the component containing $x_i.$ Note that each $G_i$ is disjoint from the set $A.$ 

Let $Q$ and $W$  be subgraphs of $G$ or subsets of $V(G).$ We say
that $Q$ and $W$ are {\it adjacent} if there exists an edge with one endpoint in $Q$ and the other endpoint in $W;$ otherwise $Q$ and $W$ are {\it nonadjacent}. Next we prove three claims.

Claim 1. $V(G_i)\cap V(G_j)=\phi,$  $G_i$ and $G_j$ are nonadjacent for $1\le i<j\le k+1.$

In the sequel, for notational simplicity, a vertex $v$ may also mean the set $\{v\}.$
We will use the fact that if $T$ is a minimum vertex cut of $G,$ then every vertex in $T$ has a neighbor in every component of $G-T.$
Clearly, $G$ has connectivity $k.$ Since $A\setminus x_j$ is a minimum vertex cut of $G,$ the subgraph $G[x_i\cup V(G_j)]$ is connected and it is
contained in the component of $H_i$ containing $x_i.$ By the definition of $G_i,$ we deduce that $(x_i\cup V(G_j))\cap V(G_i)=\phi,$ implying
$V(G_i)\cap V(G_j)=\phi.$

To show the second conclusion, just note that any vertex in $G_i$ and any vertex in $G_j$ lie in different components of the graph
$G-(A\setminus x_i).$

Claim 2. $A\cup(\mathop{\cup}\limits_{i=1}^{k+1} V(G_i))=V(G).$

To the contrary, suppose that $F\triangleq V(G)\setminus \{A\cup(\mathop{\cup}\limits_{i=1}^{k+1} V(G_i))\}$ is not empty. Let $F_1, F_2, \ldots, F_s$
be the components of $G[F].$

Recall that by definition,  for $1\le i\le k+1,$ $G_i$ denotes the union of all the components of $G-(A\setminus x_i)$
except the component $R_i$ that contains $x_i.$ Hence, for every $p$ with $1\le p\le s,$ $F_p$ is a  subgraph
of $R_i,$ implying that $G_i$ is nonadjacent to $F_p.$ Note that 
$$R_i=G[x_i\cup F\cup (\mathop{\cup}\limits_{j\ne i}V(G_j))].$$
Since $R_i$ is connected, $x_i$ is adjacent to every component of $G_j$ with $j\ne i$
and $x_i$ is adjacent to each $F_p$ for $1\le p\le s.$ Thus, every $F_p$ is adjacent to every vertex in $A.$

We choose one vertex $y_i$ from $G_i$ for each $1\le i\le k.$ Then $B\triangleq\{y_1,y_2,\ldots,y_k\}$ is an independent set of $G.$
We assert that every component of $(\mathop{\cup}\limits_{i=1}^{k+1} G_i)-B$ is adjacent to $A,$ since otherwise $G$ would have a cut-vertex.
It follows that $G-B$ is connected, contradicting condition (3). This shows that $F$ is empty and claim 2 is proved.

Claim 3. $|G_i|=1$ for every $1\le i\le k+1.$

To the contrary, suppose some $G_i$ has order at least $2.$ Without loss of generality, suppose $|G_{k}|\ge 2.$
Let $z_j$ be a neighbor of $x_{k+1}$ in $G_j$ for $j=1,\ldots,k-1.$ Since $x_{k+1}$ is adjacent to $G_k,$ $x_{k+1}$ has a neighbor $w\in G_k.$
The condition $|G_k|\ge 2$ ensures that $G_k$ has a vertex $z_k$ distinct from $w.$  Denote $C=\{z_1,z_2,\ldots,z_{k}\}.$ Then $C$ is an independent set.
We assert that every component of $(G_1\cup G_2\cup\cdots\cup G_{k})-C$ is adjacent to $A\setminus x_{k+1},$ since otherwise
some $z_j$ and $x_{k+1}$ would form a vertex cut of $G,$ contradicting the condition that $G$ is $k$-connected and $k\ge 3.$ Also,
every component of $G_{k+1}$ is adjacent to every vertex in  $A\setminus x_{k+1}.$ It follows that the graph $G-(C\cup x_{k+1})$ is connected.
But $x_{k+1}$ is adjacent to $w,$ a vertex in $G_k-z_k.$ Hence $G-C$ is connected, contradicting condition (3). This shows that each $G_i$ consists of
one vertex.

Combining the information in the above three claims, we deduce that $|G|=2k+2$ and the neighborhood of $x_i$ is $\{G_1,G_2,\ldots,G_{k+1}\}\setminus \{G_i\}$
for $1\le i\le k+1.$ It follows that $G=K_{k+1,k+1}-PM.$

Conversely, it is easy to verify that the graph $K_{k+1,k+1}-PM$ indeed satisfies the three conditions in Theorem 1. This completes the proof.
\hfill $\Box$

 Mr. Feng Liu [3] asked whether the edge version of Theorem 1 holds. The following result shows that the answer is negative.

 {\bf Corollary 2.} {\it Let $k\ge 3$ be an integer. If a graph $G$ is $k$-edge-connected with matching number greater than $k,$
 then $G$ contains a matching $M$ of cardinality $k$ such that $G-M$ is connected.}

{\bf Proof.} To the contrary, suppose that for any matching $M$ of cardinality $k,$ $G-M$ is disconnected. Consider the line graph of $G,$
denoted by $H\triangleq L(G).$ Since $G$ is $k$-edge-connected, we deduce that [5, p.283] $H$ is $k$-connected. An independent set of vertices in $H$
 corresponds to a matching in $G.$ Applying Theorem 1 to $H$  we have $H=K_{k+1,k+1}-PM,$ where we use the equality sign for graphs to mean isomorphism.
 It is known ([1] or [5, p.282]) that any line graph of a simple graph cannot have the claw as an induced subgraph.
 But for $k\ge 3,$ $K_{k+1,k+1}-PM$ contains an induced claw (many in fact). This contradiction shows that $G$ contains a matching $M$ of cardinality $k$ such that $G-M$ is connected. \hfill$\Box$

{\bf Remark.} As for the case $k=2$ of Corollary 2, using the ideas in the above proof and using the fact mentioned at the beginning of this paper,
 we see that cycles are the only $2$-edge-connected graphs in which any two nonadjacent edges form a separating set.

Finally we consider replacing independent vertices in Theorem 1 by peripheral vertices.
The {\it eccentricity} of a vertex $v$ in a graph $G,$ denoted by $e(v),$ is the distance to a vertex farthest from $v.$
A vertex  $v$ is a {\it peripheral vertex} of $G$ if $e(v)$ is equal to the diameter of $G.$  The {\it periphery} of $G$ is the set
of all peripheral vertices. We pose the following

{\bf Conjecture 3.} {\it Let $k\ge 2$ be an integer. If $G$ is a $k$-connected graph whose periphery has cardinality at least $k,$
then $G$ contains a set $S$ of $k$ peripheral vertices such that $G-S$ is connected. }

{\bf Observation 4.} {\it The case $k=2$ of Conjecture 3 is true.}

{\bf Proof.} To the contrary, suppose that any two peripheral vertices form a vertex cut of $G.$
Denote by $d(u,v)$ the distance between two vertices $u,\,v$ and let the diameter of $G$ be $d.$ We have $d\ge 2.$
Choose vertices $x,\,y$ such that $d(x,y)=d.$ Let $P$ be a shortest $(x,y)$-path, and let $y^{\prime}$ be the neighbor of $y$ on $P.$
Let $H$ be a component of $G-\{x,y\}$ that does not contain the path $P-\{x,y\}.$

Since $G$ is $2$-connected, both $x$ and $y$ have a neighbor in $H.$ Let $x^{\prime}$ be a neighbor of $x$ in $H.$ Then $d(x^{\prime}, y)\ge d-1.$
Since every $(x^{\prime}, y^{\prime})$-path contains either $x$ or $y,$ we deduce that $d(x^{\prime}, y^{\prime})=d.$ Thus $x^{\prime}$ is also
a peripheral vertex. By our assumption, $G-\{x, x^{\prime}\}$ is disconnected. Let $R$ be the component of $G-\{x, x^{\prime}\}$ containing $y.$
Clearly every component of $G-\{x, x^{\prime}\}$ other than $R$ is contained in $H.$ Let $Q$ be an arbitrary such component. We assert that every vertex in $Q$ is adjacent to $x^{\prime}.$ Let $z\in V(Q).$ Any $(z,y)$-path must contain either $x$ or $x^{\prime}.$ Since  $d(x, y)=d,$ a shortest $(z,y)$-path 
must contain $x^{\prime},$ which implies that $z$ and $x^{\prime}$ are adjacent and $z$ is a peripheral vertex, since $d(x^{\prime},y)\ge d-1.$ Choose a vertex $z_0$ from any  component of $G-\{x, x^{\prime}\}$ other than $R.$ Note that $x^{\prime}$ is adjacent to $R,$ since $\{x,x^{\prime}\}$ is a minimum
vertex cut of $G.$ Then the graph $G-\{x,z_0\}$ is connected, contradicting our assumption. \hfill$\Box$

The graph $F$ in Figure 1 shows that the connectivity condition in Conjecture 3 cannot be dropped.  $F$ has diameter $4$ and periphery
$\{ v_1,v_2,v_3,v_4,v_5,v_6\}.$ With $k=5,$ any $5$ peripheral vertices of $F$ form a vertex cut.
\vskip 3mm
\par
 \centerline{\includegraphics[width=2.7in]{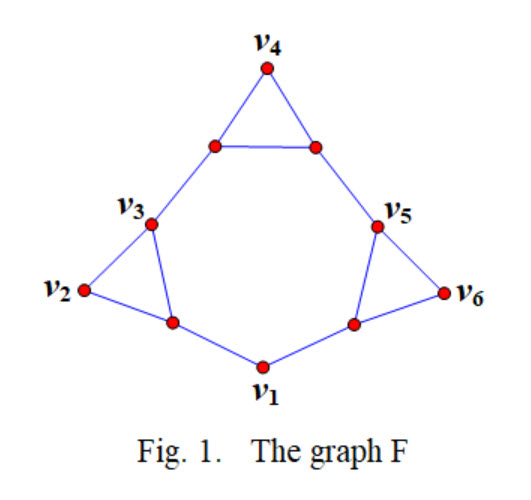}}
\par

{\bf Acknowledgement.} This research  was supported by the NSFC grant 12271170 and Science and Technology Commission of Shanghai Municipality
 grant 22DZ2229014.


\begin{thebibliography}{99}
\bibitem{1} L.W. Beineke, Derived graphs and digraphs, in Beitr\"{a}ge zur Graphentheorie, Teubner, 1968, 17-33.
\bibitem{2} J.A. Bondy and U.S.R. Murty, Graph Theory, GTM 244, Springer, 2008.
\bibitem{3} F. Liu, Private communication, October 2022.
\bibitem{4} L. Lov\'{a}sz, Combinatorial Problems and Exercises. Second Edition. North-Holland Publishing Co., Amsterdam, 1993.
\bibitem{5} D.B. West, Introduction to Graph Theory, Prentice Hall, Inc., 1996.
\end{thebibliography}
\end{document}